\documentclass[onecolumn
]{mpi2015-cscpreprint}

\usepackage[american]{babel}
\usepackage{mymacros}
\usepackage{amsbsy}

\usepackage[utf8]{inputenc}
\usepackage{pgfplots}
\DeclareUnicodeCharacter{2212}{−}
\usepgfplotslibrary{groupplots,dateplot}
\usetikzlibrary{patterns,shapes.arrows}
\pgfplotsset{compat=newest}

\usepackage{amssymb}
\usepackage{amsthm}

\usepackage[normalem]{ulem}

\newtheorem{theorem}{Theorem}[section]
\newtheorem{remark}[theorem]{Remark}

\newcommand{\bOmega}{\boldsymbol{\Omega}}
\renewcommand{\bOmega}{\Omega}
\renewcommand{\bO}{\bOmega}
\newcommand{\bGamma}{\Gamma}

\def\bOqp#1#2{\bO\,\left(#1\otimes #2\right)}

\def\JG{\partial \bGamma}
\def\JGq{\partial \bGamma_q}
\def\Mqt{M_{q(t)}}

\def\trp{\top}
\def\tt{^\top}

\def\qmf{\texttt{qmf}}
\def\dmdc{\texttt{dmdc}}

\renewcommand{\bV}{V}
\renewcommand{\bq}{q}
\renewcommand{\bx}{x}

\renewcommand{\bA}{A}
\renewcommand{\bH}{H}
\renewcommand{\bB}{B}
\newcommand{\tbA}{\widetilde{\bA}}
\newcommand{\tbB}{\widetilde{\bB}}
\newcommand{\hbA}{\widehat{\bA}}
\newcommand{\hbB}{\widehat{\bB}}

\usepackage[colorinlistoftodos,prependcaption,textsize=small]{todonotes}

\hypersetup{%
    pdftitle={A quadratic decoder approach to nonintrusive reduced-order modeling of nonlinear dynamical systems},
  pdfauthor={Peter Benner, Pawan Goyal, Jan Heiland, Igor Pontes Duff},
    pdfkeywords={operator inference, nonintrusive model reduction,
    incompressible Navier-Stokes equation}
    }

\begin{document}
  
%%%%%%%%%%%%%%%%%%%%%%%%%%%%%%%%%%%%%%%%%%%%%%%%%%%%%%%%%%%%%%%%%%%%%%%%%%%%%%%%
% PAPER INFORMATION.                                                           %
%%%%%%%%%%%%%%%%%%%%%%%%%%%%%%%%%%%%%%%%%%%%%%%%%%%%%%%%%%%%%%%%%%%%%%%%%%%%%%%%

\title{A quadratic decoder approach to nonintrusive reduced-order modeling of nonlinear dynamical systems}

  % \thanks{The first and third authors were supported by the German Research Foundation (DFG) priority program 1897: "Calm, Smooth and Smart - Novel Approaches for Influencing Vibrations by Means of Deliberately Introduced Dissipation" and the German Research Foundation (DFG) research training group 2297 "MathCoRe", Magdeburg.
% }}

\author[1]{Peter~Benner}
\affil[1]{Max Planck Institute for Dynamics of Complex Technical Systems,
  Magdeburg, Germany.\newline
Faculty of Mathematics, Otto von Guericke University Magdeburg, Germany.\authorcr
  Email: \texttt{\href{mailto:benner@mpi-magdeburg.mpg.de}{benner@mpi-magdeburg.mpg.de}},
  ORCID: \texttt{\href{https://orcid.org/0000-0003-3362-4103}%
    {0000-0003-3362-4103}}}

\author[2]{Pawan~Goyal}
\affil[2]{Max Planck Institute for Dynamics of Complex Technical Systems,
  Magdeburg, Germany.\authorcr
  Email: \texttt{\href{mailto:goyalp@mpi-magdeburg.mpg.de}{goyalp@mpi-magdeburg.mpg.de}},
  ORCID: \texttt{\href{https://orcid.org/0000-0003-3072-7780}%
    {0000-0003-3072-7780}}}

\author[3,$\ast$]{Jan~Heiland}
\affil[3]{Max Planck Institute for Dynamics of Complex Technical Systems,
  Magdeburg, Germany. \newline
Faculty of Mathematics, Otto von Guericke University Magdeburg, Germany.
  \authorcr
  Email: \texttt{\href{mailto:heiland@mpi-magdeburg.mpg.de}{heiland@mpi-magdeburg.mpg.de}},
  ORCID: \texttt{\href{https://orcid.org/0000-0003-0228-8522}%
    {0000-0003-0228-8522}}}

\author[4]{Igor~Pontes~Duff}
\affil[4]{Max Planck Institute for Dynamics of Complex Technical Systems,
  Magdeburg, Germany.\authorcr
  Email: \texttt{\href{mailto:pontes@mpi-magdeburg.mpg.de}{pontes@mpi-magdeburg.mpg.de}},
  ORCID: \texttt{\href{https://orcid.org/0000-0001-6433-6142}%
    {0000-0003-3072-7780}}}
\affil[$\ast$]{corresponding author}

\shorttitle{Operator inference on quadratic manifolds} 
%\shortauthor{Benner et al.}
\shortauthor{P. Benner, P. Goyal, J. Heiland, I. Pontes Duff}
  
\keywords{Computational fluid dynamics, scientific machine learning, Navier-Stokes equations, operator inference}

\msc{37N10, 68T05, 76D05, 65F22, 93A15, 93C10}
  
\abstract{%
  Linear projection schemes like Proper Orthogonal Decomposition can efficiently
  reduce the dimensions of dynamical systems but are naturally limited, e.g.,
  for convection-dominated problems. Nonlinear approaches have shown to
  outperform linear methods in terms of dimension reduction versus accuracy
  but,  typically, come with a large computational overhead. In this work, we consider a quadratic reduction scheme which induces nonlinear structures that are well accessible to tensorized linear algebra routines. We discuss that nonintrusive approaches can be used to simultaneously reduce the complexity in the equations and propose an operator inference formulation that respects dynamics on nonlinear manifolds.
}

\novelty{Formulation and implementation of \emph{operator inference} on a
nonlinear manifold. In particular, we

\begin{itemize}
	\item formally state the \emph{operator inference} problem on a nonlinear
    manifold including
  \item a state-dependent mass matrix that accounts for the derivative
    of the nonlinear coordinate transformation and
	\item provide a numerical example with a comparison to dynamic-mode decomposition. 
\end{itemize}
}

\maketitle

%%%%%%%%%%%%%%%%%%%%%%%%%%%%%%%%%%%%%%%%%%%%%%%%%%%%%%%%%%%%%%%%%%%%%%%%%%%%%%%%
% PAPER CONTENT.                                                               %
%%%%%%%%%%%%%%%%%%%%%%%%%%%%%%%%%%%%%%%%%%%%%%%%%%%%%%%%%%%%%%%%%%%%%%%%%%%%%%%%
  
\section{Introduction}
\label{sec:intro}

Both the idea and the promise of model order reduction are the identification and employment of coordinate systems that can best encode states and behaviors with the least amount of degrees of freedom. 
Linear projection schemes like the proper orthogonal decomposition (POD) have proven their universality and efficiency but come with limits in terms of low-dimensionality that is related to the notion of the \emph{Kolmogorov n-width}; cp. \cite{OhlR16}.
In the very low-dimensional regime, techniques that use nonlinear
relations between the actual and the reduced coordinates seem to outperform POD; see, e.g., \cite{HeiBB22, LeeC20, morGoyB21d}.
While nonlinearities may lead to fewer dimensional coordinates, their inclusion
in simulation schemes comes with extra computational efforts.% ; hence, these nonlinear techniques have not yet led to a countable efficiency gain.

A potential general approach to soften the computational disadvantages over
linear methods lies in consideration of tensor-structured nonlinearities.
These can well encode or approximate general types of nonlinear relations and
are structured to be efficiently treated by numerical linear algebra tools.
As an example and as the basic approach in this paper, we consider the quadratic relation that decodes some $r$-dimensional coordinates $q$ to approximate some $n$-dimensional coordinates $\bx$ as
\begin{equation}\label{eq:quad-decoding}
	\bx(t)  \approx \tx(t) =\bGamma(\bq(t)) : = \bV \bq(t) + \frac{1}{2}\bOmega
	\left(\bq(t) \otimes \bq(t)\right),
\end{equation}
with the matricized tensor $\bOmega\in \mathbb R^{n\times k^2}$ and where
$\otimes$ stands for the \emph{Kronecker product}. 
% to be the matrization of a symmetric tensor, where symmetry means that $\bO\,\left(\bq_1(t) \otimes \bq_2(t)\right)=\bO\,\left(\bq_2(t) \otimes \bq_1(t)\right)$. 

%\igor{\sout{As an example application, we consider the quadratic model systems of the form}} 
In this work, we focus on data-driven reduced order modeling of quadratic model systems of the form 
\begin{equation}\label{eq:quad_dyn}
  \dot{\bx}(t) = \bA \bx(t) +\bH \left(\bx (t) \otimes \bx(t)\right) + \bB 
	, \quad \bx(0)=x_0,
\end{equation}  
with the (full order) state $\bx(t) \in \R^{n}$. 
% We will assume w.l.g. that $\bH$ is a symmetric tensor, i.e., $\bH(\bx_1 \otimes \bx_2) = \bH(\bx_2 \otimes \bx_1)$.
%
As a benchmark and to introduce the basic concepts, we briefly discuss the
method of \emph{Proper Orthogonal Decomposition} (POD): Given
a matrix $X \in \mathbb R^{n\times k}$ of $k$ solution snapshots
\begin{equation}\label{eq:snapshot-matrix}
	X=
	\begin{bmatrix}
		x(t_1) ,~x(t_2) ,~\hdots,~x(t_k)
	\end{bmatrix}
  ,
\end{equation}
then POD identifies a matrix $V\in \mathbb R^{n\times r}$ with $r\leq k$, with
$V^\trp V = I_r$, so that, measured in the \emph{Frobenius norm} $\|\cdot \|_F$,
the projection error
% \todo{I have added this orthogonality condition. } so that
\(
\|X - VV\tt X\|_F,
\)
is minimal over all possible matrices $V\in \mathbb R^{n\times r}$. This choice
means that the coordinates 
\begin{equation}\label{eq:pod-coordinates}
	q(t) = V^\trp x(t)
\end{equation}
encode the state $x$ in an $r$-dimensional space optimally on average with
respect to the linear decoding
\begin{equation*}
	\tx(t) = Vq(t)
\end{equation*}
and with respect to the data collected in $X$. If, in \eqref{eq:quad_dyn},
the state $x$ is replaced by the parametrization $\tx = Vq$ and the equations
are projected onto the span of $V^\trp$, then the POD reduced order model is
obtained as

\begin{equation*}%\label{eq:quad_dyn_pod}
	 \dot{q}(t) =V^\trp  \bA V {q}(t) +V^\trp \bH \left(V {q} (t)
	\otimes V {q}(t)\right) + V^\trp \bB 
	, \quad q(0)=V^\trp x_0.
\end{equation*}  

The motivations for this work are: 
\begin{itemize}
	\item 
	to follow the POD approach to identify low-dimensional coordinates from 
	data
	\item 
	with a quadratic decoding
	\item 
	to be used for \emph{nonintrusive model reduction} by means of
	\emph{operator inference}
	\item 
	with a nonlinear projection of the state equations.
\end{itemize}

The idea of using quadratic decoding has been proposed in \cite{JaiTRR16} for
mechanical systems, where the reduced-order coordinates can be derived from first
principles. In \cite{BarF22}, a data-driven approach for finding the coordinates
has been investigated and applied in an \emph{intrusive} way as an enhancement
of POD. As we will illustrate here, this \emph{intrusive} use of nonlinear
decodings induces higher-order nonlinearities. Thus, as we will argue, the
combination with \emph{operator inference} (see, e.g., \cite{morPehW16}), seems
to make the best use of nonlinear decoding and reduced-order modelling.
This idea
has been followed in \cite{GeeWW22} though in a \emph{one-sided} fashion,
meaning that the nonlinear embedding is used for decoding the 
coordinates but not respected in the projection of the equations. % Thus, the
% approach presented in \cite{GeeWW22} basically infers a quadratic reduced order
% system with an improved embedding in the original coordinate space. 
Thus, one can also interpret the methodology in \cite{GeeWW22} as obtaining a
classical operator inference \cite{morPehW16} with a quadratic nonlinearity and,
additinally, a quadratic correction term for decoding the original coordinates from POD coordinates.

The paper is structured as follows. In \Cref{sec:proj-mor-nonl-decod}, we lay
out how an intrusive model reduction based on nonlinear decodings would look
like. 
% We will argue that the intrusive approach increases the degree of nonlinearities
% and advocate for 
The \emph{nonintrusive} approach of inferring the operators for a system
approximation of a lower degree is explained in \Cref{sec:qmf-quad-opinf} with
implementation details provided in \Cref{sec:qmf-implementation}. 
In \Cref{sec:num-exa}, we present the setup and numerical results for an example
application of two-dimensional incompressible flows before concluding the paper
with a summarizing discussion in \Cref{sec:conclusion}.

\section{Intrusive Model Reduction with Nonlinear Decoding}
\label{sec:proj-mor-nonl-decod}

We will refer to a \emph{(quadratic) manifold} of dimension $r$, when talking about a domain
$\mathcal Q \subset \mathbb R^{r}$ that is embedded into $\mathbb R^{n}$ by a
(quadratic)
map $\Gamma\colon \mathcal Q \mapsto \mathbb R^{n}$ that is differentiable and
onto with a differentiable left inverse. 
% As for the matrix calculus, we use $\|\cdot\|_F$ \todo{PG: Since $\|\cdot\|_F$ already appears in the introduction, would it make sense to move there? Also, Kronecker-product symbol!} to denote the \emph{Frobenius}-norm (which reduces to the $2$-norm for a vector), and we use $\otimes$ to denote the \emph{Kronecker}-product.

Assume that a low-dimensional manifold $\mathcal Q \in \mathbb R^{r}$ is given
that well encodes the state of \eqref{eq:quad_dyn} with the embedding $\Gamma
\colon \mathbb R^{r}\to \mathbb R^{n}$. Then, replacing $x(t)$ by $\tilde x(t) =
\Gamma (q(t))$ in \eqref{eq:quad_dyn} and respecting the chain rule of
differentiation, we obtain the dynamical system in the reduced
coordinates as
\begin{equation}\label{eq:IncludManifold}
	 \partial \bGamma(\bq(t)) \dot q(t) = \bA \bGamma(\bq(t))
	+ \bH(\bGamma(\bq(t))\otimes \bGamma(\bq(t))) +\bB.
\end{equation}

Since $\Gamma$ is onto, it holds that its Jacobian $\JG(q(t))=:\JGq$ has full column rank, so that %, with $\bM$ being symmetric positive definite, 
a premultiplication by $\JG(q(t))^\trp$ turns
\eqref{eq:IncludManifold} into a regular ODE: %\PG{I tried top instead of trp to get transpose and i feel top comes better, i.e, $\JGq^\top$ (with top) and $\JGq^\trp$ (with trp). what do you think about replacing trp with top, but I do not have a strong opinion?}
\begin{equation}\label{eq:QM_rom}
	\JGq^\trp  \JGq \, \dot \bq(t) = \JGq^\trp\, \bA\bGamma(\bq(t)) +  \JGq^\trp\,
	\bH\, \left(\bGamma(\bq (t)) \otimes \bGamma(\bq(t))\right) +  \JGq^\trp\,  \bB
\end{equation}

\begin{remark}
	As for the choice of the projection, we would like to make the following
  relevant remarks:
	\begin{itemize}
		\item Notice that the projected system \eqref{eq:QM_rom} possesses a  matrix $\JGq^\trp \JGq$ which is full rank. Hence, the system \eqref{eq:QM_rom} is a regular ODE. This justifies the use of a $q(t)$ dependent projection. It is worth noticing that if the reduced model was obtained by enforcing Petrov-Galerkin with the constant matrix $W \in \mathbb R^{n\times r}$, the projected mass matrix would have the form $W^\trp \JG(q(t) )$ and might not be necessarily invertible for all $t$.    
		%
		%\igor{\sout{The use of a $q(t)$ dependent projection is needed, because a constant matrix $V\in \mathbb R^{n\times r}$ so that $V^\trp \bM \JG(q(t) )$ is invertible for all $t$, might not necessarily exist in general.}} \todo{I have changed this part. Let me know what you think}
		%
    \item Doing so, as in the standard POD projection approach, we also make sure that $\dot q(t)$ is the minimizer
      to
      \begin{equation*}
      \min_{\zeta \in \mathbb R^{r}}
	 \| \partial \bGamma(\bq(t)) \zeta - \bA \bGamma(\bq(t))
	- \bH(\bGamma(\bq(t))\otimes \bGamma(\bq(t))) -\bB \|_F.
      \end{equation*}
		\item  The choice of $\JG(q(t) )^\trp$ as a projection matrix can also be motivated 
		by $\JG(q(t) )$ describing a basis for the tangent space to the manifold
		$\Gamma(\mathcal Q)$ at $\Gamma(q(t))$.
	\end{itemize}
\end{remark}

% \paragraph{Model reduction with quadratic manifold maps}

If the map $\Gamma$ is quadratic in form of \eqref{eq:quad-decoding}, its
Jacobian at point $q(t) $ is given as the linear map
\begin{equation}\label{eq:jacobian-gamma}
	\tq \mapsto V\tq + \frac 12 ( \bOqp {q(t)}{\tq} + \bOqp {\tq}{q(t)})
\end{equation}
that can be realized as a matrix as in the formulas provided in \cite[Sec. 7.2]{BehBH17}. 

From \eqref{eq:jacobian-gamma}, it can be seen that $\JG(q(t))$ itself is linear
in $q(t)$, so that premultiplication by $\JGq^\trp$ will increase the degrees of
the polynomials in the states $q(t)$, e.g., the projected constant term will be
realized as
\begin{equation*}
	\JGq^\trp \bB = \tbB_1 +  \tbB_2 q(t) .
\end{equation*}
%\todo{PG: Maybe we discussed, it would be good to have hat for reduced quantities? \igor{IP: I guess wel could have a hat for inferred operators and a tilde for reduced ones. What do you think? }}

Overall, the nonlinearly reduced model using the quadratic manifold projection
that approximates \eqref{eq:quad_dyn} will have the form
\begin{equation}\label{eq:NL_rom}
	% \boxed{	
	\Mqt \dot\bq(t) =  \tbA_1\bq(t) + \tbA_2 \otimes^2 \bq(t) +  \tbA_3\otimes^3\bq(t) +  \tbA_4\otimes^4\bq(t)+ \tbA_5 \otimes^5\bq(t) +  \tbB_1,
	% }
\end{equation}
where $\otimes^m q(t):=\bigotimes_{i=1}^m q(t)$, with the $q(t)$ dependent mass matrix $\Mqt:=\JGq^\trp\JGq$, and with coefficients that, in an intrusive MOR scheme, have to be derived from the original coefficients $\bA$, $\bH$, and
$\bB$ and their products with $\JG(q(t) )^\trp$. 
%\igor{Notice that the quantities $\tbA_1, \tbA_2, \tbA_3, \tbA_4, \tbA_5,$ and $\tbB_1$ are obtained intrusively, i.e., it is required the knowledge of the full order model matrcies  $\bA$, $\bH$, and, $\bB$.   }

We will use, however, the approach of \emph{operator inference} to extract
lower-order coefficients that best approximate \eqref{eq:NL_rom} with respect to
given trajectory data.

Apart from being purely data-driven, i.e., no knowledge of the original
coefficients is required, the use of \emph{operator inference} for a model with
a lower polynomial degree comes with the
conceptual advantage that higher orders are possibly compensated by them.
Thus, only coefficients of moderate size can be inferred, stored, and
later evaluated, and possible numerical instabilities that come with large
exponents are avoided. For this reason, in this work, we will limit ourselves to surrogate models of the form

\begin{equation}\label{eq:Sor_rom}
	% \boxed{	
	\Mqt \dot\bq(t) = \tbA_1\bq(t) + \tbA_2 \otimes^2 \bq(t) + \tbB_1.
	% }
\end{equation}

\section{Operator Inference for Low-order
Approximation}\label{sec:qmf-quad-opinf}

In this section, we propose a non-intrusive method enabling to infer surrogate
models for \eqref{eq:quad_dyn} using a quadratic encoding as in
\eqref{eq:quad-decoding}. Our main goal is to infer reduced-order models of the form \eqref{eq:Sor_rom}. To this aim, we will require two steps. The first step is to construct the quadratic encoding using snapshots of the state. Then, based on this quadratic encoding, we will infer this reduced system operators by means of a tailored operator inference approach.

We assume that we are given a matrix $X \in \mathbb R^{n\times k}$ of $k$ solution snapshots of a dynamical system
with state $x(t)\in \mathbb R^{n}$ on a grid $t_1 < t_2 < \dotsm <
t_k$ as in \eqref{eq:snapshot-matrix} % \todo{PG: I am not sure if we need this assumption  of data collected at equidistant?}
and a POD based encoding and decoding
\begin{equation*}
  q(t) = V^\trp x(t) \quad \text{and} \quad \tx(t) = Vq(t)
\end{equation*}
as in \eqref{eq:pod-coordinates} with $V\in \mathbb R^{n\times r}$. 

As in \cite{BarF22} and \cite{GeeWW22}, we define the quadratic decoding as a
correction of the POD decoding, which can be formulated as $\Omega \in \mathbb
R^{n\times k^2}$ being a minimizer for 
\begin{equation}\label{eq:inferring-Omega}
  \min_{\Omega \in \mathbb R^{n\times r^2}}
  \|X - VQ - \frac 12 \Omega \, ( Q\otimes_c Q) \|_F^2,
\end{equation}% \todo{IP: Should we change $\mathcal O$ to $\Omega$ in the optimization? }
with the matrix 
$$V^\trp X=:Q:= \begin{bmatrix}q_1 ,~ q_2 ,\hdots, ~q_r
\end{bmatrix}\in \mathbb R^{r\times k}$$
and where we define 
\begin{equation*}
  Q \otimes_c Q = \begin{bmatrix}
    q_1 \otimes q_1,~ q_2 \otimes q_2,\hdots,~ q_k \otimes q_k
  \end{bmatrix} \in \mathbb R^{r^2 \times k}.
\end{equation*}

With $\Omega$ at hand, we realize the gradient $\JG{q_j}$ as in
\eqref{eq:jacobian-gamma} at all snapshots $j=1,\dotsc,k$, and, having approximated
the time derivatives $\dot q_j$ of the snapshots $q_j$, we can set up the matrix 
\begin{equation*}
  \dot Q_M := \begin{bmatrix}
    M_{q_1}\dot q_1 ,~M_{q_2}\dot q_2 ,\hdots ,~ M_{q_r}\dot q_r
  \end{bmatrix},
\end{equation*}
with $M_{q_k} = \partial \bGamma_{q_k}\partial \bGamma_{q_k}$,
and state the operator inference problem as
\begin{equation}\label{eq:qmf_quad_opinf_problem}
  \min_{\widehat{\mathcal A}_1 \in \mathbb R^{r\times r},\, 
    \widehat{\mathcal A}_2 \in \mathbb R^{r \times r^2},\,
    \widehat{\mathcal  B}_1 \in \mathbb R^{r\times 1}
  }
  \| \dot Q_M - \widehat {\mathcal A}_1 Q - \widehat{\mathcal A}_2\,( Q\otimes_c Q) - \widehat{ \mathcal
  B}_1\|_F^2.
\end{equation}

Given (approximately) minimizing solutions $\hbA_1$, $ \hbA_2$, and $ \hbB_1$ to
\eqref{eq:qmf_quad_opinf_problem}, a lower-order approximation to
\eqref{eq:NL_rom}, and via $x(t) \approx \tx(t)=Vq(t) + \frac 12 \Omega \,
q(t)\otimes q(t)$, a reduced-order model of the dynamics of \eqref{eq:quad_dyn},
is defined through
\begin{equation*}
  M_{q(t)}\dot q(t) =  \hbA_1 q(t) + \hbA_2 \, \left(q(t)\otimes q(t)\right) + \hbB_1, \quad
  q(0)=V^\trp x(0).
\end{equation*}

\section{Implementation Details}
\label{sec:qmf-implementation}
Regarding the numerical implementation of the \emph{operator inference} approach
of \Cref{sec:qmf-quad-opinf}, we address two issues.

Firstly, for a vector $q\in \mathbb R^{r}$, the Kronecker product $q\otimes q \in \mathbb
R^{r^2}$ has only $\frac{(r+1)r}{2}$ different entries. Accordingly, the function
$\Omega\, (q\otimes q)$, with $\Omega \in \mathbb R^{r\times r^2}$, can be
equivalently realized by a matrix $\tilde \Omega \in \mathbb R^{r\times
	\frac{(r+1)r}{2}}$ and a corresponding removal of the redundant entries of
$q\otimes q$. 
We make use of this compression by directly inferring the reduced realizations
of $\Omega$ and $A_2$ in \eqref{eq:inferring-Omega} and~\eqref{eq:qmf_quad_opinf_problem}.

Secondly, it has been commonly observed that an overfitting of the coefficients
in \eqref{eq:inferring-Omega} and \eqref{eq:qmf_quad_opinf_problem} has
adversary effects. The best results are obtained by trading in some accuracy for a
smaller norm of the solutions to the optimization problems. When using the SVD
to compute the coefficients, this tradeoff can be achieved by truncating the SVD
at suitably chosen threshold values; see \cite{BenGHP22} for a heuristic
strategy for choosing these parameters based on \emph{L-curves}. 
We further note that this truncation has a similar effect as a regularization of the optimization problem as done in
\cite{GeeWW22}.

\section{Numerical Example}\label{sec:num-exa}
\begin{figure}[tb]
	\includegraphics[width=\textwidth]{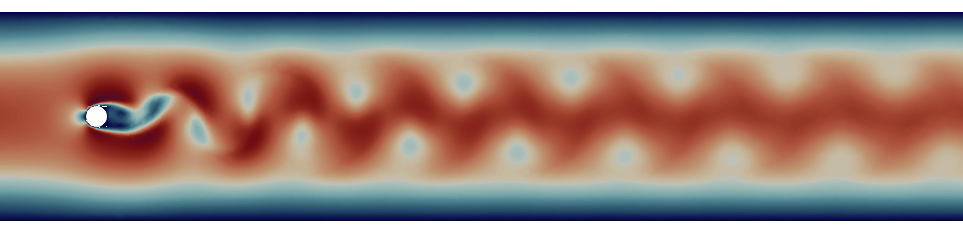}
	\caption{Snapshot of the velocity magnitude of the considered setup in the
		fully developed transient regime of $\mathsf{Re}=30$.}
	\label{fig:sol-snapshot-re30}
\end{figure}

We consider the flow past a cylinder in 2D at Reynolds number $\mathsf{Re}=30$ 
% \todo{PG: In the figure Re = 60 but here 30. Can you please check?} 
in
the startup phase from the associated steady-state \emph{Stokes} solution and in the fully
developed transient regime. 
The geometrical setup and a snapshot of the developed flow is presented in
\Cref{fig:sol-snapshot-re30}; for a detailed description see
\cite[Sec. 5]{HeiBB22}.

In the spatial domain, we use a \emph{Finite Element} discretization that approximates the flow velocity with $41,682$ degrees of freedom.
We consider a time frame from $t=0$ to $t=12$ and take $2000$ snapshots of the
velocity solution $v$ at equispaced time instances that cover the whole time
range.

\begin{figure}[tb]
	\input{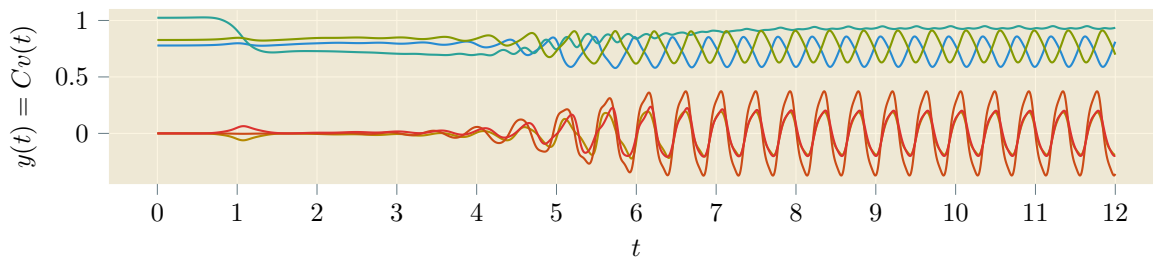}
	\caption{Output of the full order simulation.}
	\label{fig:output-fom}
\end{figure}

The approximation will target the actual system states $v(t)$, and we will report
the numerical realization of the approximation error
\begin{equation*}
	% e(t_0, t_e) = \int_{t_0}^{t_e} 
	\|v(t) - \tilde v(t)\|_{L^2(\mathcal
		D;\mathbb R^{2})},
	%\operatorname{d}t,
\end{equation*}
where $(t_0, t_e)$ defines the time interval, $v(t)$ is described by the
snapshots, $\tilde v(t)$ is the reconstruction from the reduced order
model, and  $\mathcal D \subset \mathbb R^{2}$ denotes the computational
domain. For illustration purposes, we report an output $y(t)=Cv(t)\in \mathbb
R^{6}$ that is obtained by measuring the locally spatially averaged values of the two
components of the velocity at three \emph{sensor locations} in the cylinder
wake. The output of the full order simulation is depicted in
\Cref{fig:output-fom}, showing
the development of the flow from the initial state to the characteristic
periodic vortex-shedding regime.

We check the performance of  the following two methods:
\begin{itemize}
	\item[\qmf] --- the quadratic manifold approach as developed in
	\Cref{sec:qmf-quad-opinf} and
	\item[\dmdc] --- the \emph{Dynamic Mode Decomposition} with a bias term; see,
	e.g., \cite{morProBK16}.
\end{itemize}
Both methods are completely data-driven, and \dmdc~has shown to be a suitable
benchmark because of its good performance independent of hyperparameters other
than the dimension $r$ of the reduced order model; see
\cite[Sec. 7.2]{BenGHP22}

As none of the methods could encode the full transition of the rather
nonphysical \emph{Stokes} steady state towards the fully developed periodic
regime, we checked the performance for the initial phase for $t\in (0, 6)$
and the final phase for $t\in (6,12)$ separately with separate models. In both
cases, the first $500$ of the $1000$ snapshots that fell into the corresponding
subinterval
% \todo{But earlier it is stated that there are total 2000 snapshots?} 
were used to infer the models and the remaining snapshots served as an estimation for the prediction capabilities of the approaches.

The approximation results for $r=10$ and truncations of the SVDs that solve
\eqref{eq:inferring-Omega} and \eqref{eq:qmf_quad_opinf_problem} to keep only
the first $7$ and $10$ largest singular values, respectively, are plotted in
\Cref{fig:r10}.

In the initial phase, both methods perform similarly. While an initial transient
response is well captured, the developing transition to the periodic regime is
not reflected in the outputs, and the approximation error suddenly grows once the
subinterval that was covered by the training data is left. In the transient phase, both approaches show good
extrapolation capabilities. Overall, on this example, the \qmf~approach well competes with
established methods. Additional potentials and good adaptivity to other problem
classes are to be provided, e.g., by a systematic
optimization of the hyperparameters.
% , with the \qmf~approach performing slightly better. % \todo{PG: I  am not sure: A picky reviewer might say dmdc also works better for some time in the testing regime }
\begin{figure}
	\input{figures/sc_30_0.0-5.996_2000_rv10.tex}
	\input{figures/sc_30_6.003-12.0_2000_rv10.tex}
	\caption{Outputs and approximation errors for the \qmf~and \dmdc~model of reduced order
		$r=10$ in the initial phase (left) and the periodic regime (right). The thin
		black lines in the upper plots denote the output of the full order model; cp.
		\Cref{fig:output-fom}. The dashed vertical lines in the lower plots separate the
		regions of the training data (left of the dashed line) and the region that is
		extrapolated by the model.}
	\label{fig:r10}
\end{figure}

The code and the raw data of the presented numerical results are available as
noted in \Cref{fig:linkcodedata}.

\section{Conclusion and Discussion}\label{sec:conclusion}

The presented approach to constructing reduced-order models is purely data-driven. Future theoretical research could
address system analytical investigations
of low-dimensional manifolds that encode system states, e.g., as
in \cite{KolW20} or rigorous reasoning, as in \cite{JaiTRR16}, that
for the system under consideration, a quadratic
manifold suggests itself for efficiently encoding the state.

The presented numerical results show the potential but no clear advantage over the very robust \dmdc~method. Future work would seek algorithmic improvements like the choice of different norms for the involved optimization
problems or reliable strategies to choose hyperparameters.
In a finite-element context, e.g., the corresponding discrete Sobolev norms will provide
targeted measures for the optimizations that are consistent with the underlying
PDE model and useful for estimating truncation errors independent of other model
parameters.

%%%%%%%%%%%%%%%%%%%%%%%%%%%%%%%%%%%%%%%%%%%%%%%%%%%%%%%%%%%%%%%%%%%%%%%%%%%%%%%%

\begin{figure}[ht]
  \caption{Code and Data Availability.}\label{fig:linkcodedata}
  \begin{center}
    \fbox{
      \begin{minipage}{.9\linewidth}
        The raw data and the source code of the implementations used to compute the presented results is available from
        \begin{center}
          \href{https://doi.org/10.5281/zenodo.7126187}{\texttt{doi:10.5281/zenodo.7126187}}
        \end{center}
        under the MIT license and is authored by Jan Heiland.
      \end{minipage}
    }
  \end{center}
\end{figure}

\section{Acknowledgements}
Peter Benner, Jan Heiland, and Igor Pontes are supported by the German Research
Foundation (DFG) through the Research
Training Group 2297 “MathCoRe”, Magdeburg.

%%%%%%%%%%%%%%%%%%%%%%%%%%%%%%%%%%%%%%%%%%%%%%%%%%%%%%%%%%%%%%%%%%%%%%%%%%%%%%%%
% *** REFERENCES ***                                                           %
%%%%%%%%%%%%%%%%%%%%%%%%%%%%%%%%%%%%%%%%%%%%%%%%%%%%%%%%%%%%%%%%%%%%%%%%%%%%%%%%

\addcontentsline{toc}{section}{References}
% \bibliographystyle{abbrvurl}
% \bibliography{csc,mor,qmf}

\end{document}